\documentclass[11pt,twoside,a4paper]{article}
\usepackage{amsmath,amssymb,amsthm}
\newtheorem{Thm}{Theorem}[section]
\newtheorem{Lem}[Thm]{Lemma}
\newtheorem{Pro}[Thm]{Proposition}

\theoremstyle{definition}

\theoremstyle{remark}
\newtheorem{Rem}[Thm]{Remark}

\newcommand{\R}{\mathbb{R}}

\newcommand{\N}{\mathbb{N}}

\newcommand{\cP}{\mathcal{P}}

\newcommand{\ep}{\varepsilon}
\newcommand{\de}{\delta}
\newcommand{\De}{\Delta}
\newcommand{\ga}{\gamma}
\newcommand{\Ga}{\Gamma}
\newcommand{\al}{\alpha}

\newcommand{\si}{\sigma}

\newcommand{\la}{\lambda}

\renewcommand{\d}{\partial}

\newcommand{\dinf}{\d_{\infty}}

\newcommand{\sm}{\setminus}
\newcommand{\sub}{\subset}
\newcommand{\es}{\emptyset}

\newcommand{\ov}{\overline}
\newcommand{\wt}{\widetilde}

\newcommand{\dist}{\operatorname{dist}}

\newcommand{\hyp}{\operatorname{H}}

\newcommand{\vol}{\operatorname{vol}}
\newcommand{\iso}{\operatorname{Iso}}

\begin{document}

\title{An estimate for the volume entropy of nonpositively
curved graph manifolds}

\author{S.~Buyalo
\footnote{Partially supported by RFFI Grants 99-01-00104
and 00-15-96024.}}
\date{}
\maketitle

\begin{abstract}
Let
$M$
be a closed 3-dimensional graph manifold. We prove that
$h(g)>1$
for each geometrization
$g$
of
$M$,
where
$h(g)$
is the topological entropy of geodesic flow of
$g$.
\end{abstract}

\section{Introduction} Asymptotic geometry of a nonpositively
curved (NPC for brevity) graph manifold is a complicated mixture
of flat and hyperbolic parts which both nontrivially contribute
in general picture. Recall that a NPC metric on a closed
3-dimensional graph manifold
$M$
recovers the JSJ-decomposition of
$M$
in the following sense. There is the unique (up to isotopy)
minimal finite collection
$E$
of flat geodesically embedded tori and Klein bottles which
are pairwise disjoint, and the metric completion
of each connected component of the complement to
$E$
is a Seifert space called a block of 
$M$.
Each block
$M_v$
is fibered over a 2-orbifold
$S_v$
with negative Euler characteristic,
$\chi(S_v)<0$.
Furthermore, the metric locally splitts as
$U\times (-\ep,\ep)$
along the interior of each block, where
$U$
is a NPC surface, the splitting is compatible 
with the fibration, the fibers are closed geodesics, and
the regular fibers have one and the same length 
$l_v>0$
depending only on the block.

Since we are interested in asymptotic properties, which
are certainly the same for any finite covering of
$M$,
we may assume for simplicity that 
$M$
is orientable, the collection
$E$
consists of tori, and each block is a trivial 
$S^1$-bundle
over a compact surface
$S_v$
with boundary,
$M_v=S_v\times S^1$.
We also assume that the graph manifold structure of
$M$ 
is nontrivial, i.e.
$M$
itself is not a Seifert fibered space (though it may consist of
one block).

The flat part of the asymptotic geometry of 
$M$
studied in \cite{BS}, \cite{CK}, see also \cite{HS}. 
Roughly speaking, it can be described by very special
geodesic rays
$[0,\infty)\to M$,
which terminate in no block, skip through separating tori
$e\in E$
almost tangently spending most of the time near tori and
moreover this time rapidly increases with each step.
Though the set of such rays is a negligible part all of the
rays, it contains an important information about 
geometry of
$M$:
it was shown in \cite{CK} how this information allows to recover 
(up to scaling) the marked length spectrum of closed geodesics on
$S_v$
and the fiber length for each block
$M_v$.

Here we study the hyperbolic part of the asymptotic
geometry of
$M$,
assuming that the surface
$U$
from the local splitting
$U\times(-\ep,\ep)$
above has the constant curvature
$K=-1$.
In other words, each block fibers over a {\it hyperbolic}
orbifold (surface)
$S_v$.
A NPC metric on 
$M$
satisfying this condition is called a {\it geometrization} of
$M$
(because the metric of each block is modelled on
$\hyp^2\times\R$). Note that any geometrization of
$M$
is only 
$C^{1,1}$-smooth
being analytic along the interior of each block. It is 
known \cite{L} that 
$M$
admits a NPC metric if and only if it admits a 
geometrization. Necessary and sufficient topological
conditions for
$M$
to carry a NPC metric were found in \cite{BK}.

The relevant metric invariant which measures hyperbolicity
of a space is the volume entropy 
$h$.
Let
$\pi:X\to M$
be the universal covering,
$x_0\in X$.
Recall that
$h=h(X)$
is defined by
$$h=\lim_{R\to\infty}\frac{1}{R}\ln\vol B_R(x_0),$$
where
$B_R(x)$
is the ball in
$X$
of radius
$R$
centered at
$x$.
It is well known \cite{M} that the limit exists,
$h$
is independent of the choice of
$x_0$
and if
$M$
is NPC then
$h$
coincides with the topological entropy of geodesic flow of
$M$.
The volume entropy scales as
$l^{-1}$,
$l$
is the length. Thus the choice of a geometrization 
$g$
of
$M$
also serves as a normalization. The sectional curvatures of
$g$
satisfy
$-1\le K\le 0$.
Hence 
$h(g)\le 2$
by comparison with
$\hyp^3$.
Our main result is this.

\begin{Thm}\label{Thm:main} For any geometrization
$g$
of a graph manifold
$M$
we have
$h(g)>1$.

\end{Thm}

\begin{Rem} Though the universal covering
$X$
of
$M$
looks much more complicated than the model space
$\hyp^2\times\R$,
even the estimate
$h(g)\ge 1=h(\hyp^2\times\R)$
is not obvious and nontrivial:
$X$
contains no isometrically and geodesically embedded
$\hyp^2$, which would lead to
$h(g)\ge 1$;
on the other hand, the attempt to compare
$X\to\hyp^2\times\R$
via exponential maps identifying some tangent spaces
fails because the Jacobian of this map is
$>1$
at some points. Finally, the estimate from \cite{BW} for
the measure theoretic entropy of geodesic flow, which is always
$\le h(g)$,
gives only 
$\pi/4<1$
as a lower bound for any geometrization 
$g$
of
$M$
(if one ignores the fact that
$C^{1,1}$-smoothness
of
$g$
is not sufficient to apply this estimate).
\end{Rem}

To prove Theorem~\ref{Thm:main} we use the well known fact that
$h(g)$
coincides with the critical exponent of a Poincar\'e series
$$\cP(t)=\sum_{\ga\in\Ga}e^{-t|x_0-\ga x_0|},$$
where the fundamental group
$\Ga=\pi_1(M)$
isometrically acts 
on
$X$
as the deck transformation group. Actually, instead of
$\cP$
we use a modified
Poincar\'e series
$\cP_W$,
where summation is taken over a set
$W$
of walls in
$X$. Our proof involves three ingredients: (i) a local
estimate, which is technical and used in (ii); this estimate
is obtained in section~\ref{sect:locest}.
(ii) An accumulating procedure, which is an inductive construction 
of appropriately choosen broken geodesics in
$X$
between the base point
$x_0$
and the walls from
$W$;
the choice of these paths is the key point of the proof. The outcome
of the accumulating procedure is a generating set for
$\cP_W$
used in (iii); the procedure is described in sect.~\ref{sect:accpro}.
(iii) A self-similarity type argument. This part of the proof
uses a standard idea from self-similarity theory to show that
$\cP_W(h)$
diverges for some 
$h>1$
using the generating set obtained in (ii). This is done in 
sect.~\ref{sect:selfsim}.

\medskip
{\bf Acknowledgment.} The author is grateful to W.~Ballmann for
useful discussions of the topic of this note.

\section{Local estimate}\label{sect:locest}

Let
$F$
be the universal covering of a compact hyperbolic surface
$S$
with geodesic boundary. We identify
$F$
with a convex subset
$F\sub\hyp^2$
bounded by countably many disjoint geodesic line and fix
$o\in\hyp^2\sm F$.
Let
$w_0$
be the boundary line of
$F$
closest to
$o$, $o_0\in w_0$
be the point on
$w_0$
closest to
$o$,
so that
$|o-o_0|=\dist(o,F)=:l>0$.

Let
$A$
be the set of the boundary lines of
$F$
different from
$w_0$.
For each
$w\in A$
we denote by
$o_w\in w$
the point closest to
$o$.
Then the geodesic segment
$oo_w$
intersects
$w_0$
at some point
$t_w$,
and for
$\wt l_w=|o-o_w|$
we have
$$\wt l_w=\wt l_w'+l_w'',$$
where
$\wt l_w'=|o-t_w|$, $l_w''=|t_w-o_w|$
(all distances are taken in
$\hyp^2$).

Next, we identify
$\hyp^2$
with
$\hyp^2\times 0\sub\hyp^2\times\R$,
so that
$F$
becomes a subset of
$\hyp^2\times\R$,
and we use for it the notations introduced above. Note that
the point
$o_0$
is closest to
$o$
from the wall
$w_0\times\R$.
We take a geodesic line
$\si\sub w_0\times\R$ 
through
$o_0$
which is not horizontal, i.e.,
$\si\not=w_0\times 0$,
and take
$s_w\in\si$
with
$|s_w-o_0|=|t_w-o_0|$.
Now, we put
$l_w':=|o-s_w|$
(the distance is taken in
$\hyp^2\times\R$),
$\De_w:=\wt l_w'-l_w'$.

In other words, we replace the distance
$\wt l_w'$
between
$o$
and
$t_w$
is the hyperbolic plane
$\hyp^2$
by the distance
$l_w'$
in
$\hyp^2\times\R$,
which is shorter by comparison: the triangles
$oo_0t_w\sub\hyp^2\times 0$, $oo_0s_w\sub\hyp^2\times\R$
both have the right angles at
$o_0$, $\angle(oo_0t_w)=\frac{\pi}{2}=\angle(oo_0s_w)$,
the common side
$oo_0$
and equal sides
$|o_0-t_w|=|o_0-s_w|$.
Since
$oo_0t_w$
lies in the hyperbolic plane
$\hyp^2\times 0$
but
$oo_0s_w$
not,
we have
$\De_w>0$
except the case
$t_w=o_0=s_w$.
Now, we want to estimate 
accumulation of the 
differences
$\De_w$ 
from below. The precise statement is this.

\begin{Lem}\label{Lem:locest} Given
$l_0>0$, $\al_0\in(0,\pi/2]$, there exists
$\la_0>1$, 
which depends only on
$l_0$, $\al_0$
and the compact surface
$S$,
so that
$$\la(F,l,\al):=e^l\sum_{w\in A}e^{\De_w}e^{-\wt l_w}\ge\la_0,$$
whenever
$l=\dist(o,w_0)\ge l_0$
and the angle
$\al$
between the lines
$w_0\times 0$
and
$\si$
is at least
$\al_0$, $\al_0\le\al\le\pi/2$.
\end{Lem}

\begin{proof} By well known formula of hyperbolic geometry we have
$e^l=(\tan\frac{\psi}{4})^{-1}$,
$e^{-\wt l_w}=\tan\frac{\psi_w}{4}$,
where
$\psi$, $\psi_w$
are the angles under which
$w_0$
respectively
$w\in A$
are observed in
$\hyp^2$
from
$o$.
The boundary at infinity
$\dinf F\sub\dinf\hyp^2=S^1$
coincides with the limit set of
$\pi_1(S)$
represented in
$\iso(\hyp^2)$
as a Fuchsian group of second kind. It is well known
that the Hausdorff dimension of
$\dinf F$
(with respect to the angle metric) is 
$<1$,
in particular, the Lebesgue measure of
$\dinf F$
is zero. Thus
$\psi=\sum_{w\in A}\psi_w$.
Therefore,
$\tan\frac{\psi}{4}\le\sum_{w\in A}\tan\frac{\psi_w}{4}$
and
$$e^l\sum_{w\in A}e^{-\wt l_w}=\sum_{w\in A}\tau_w\ge 1,$$
where
$\tau_w=e^{l-\wt l_w}$.
However, this sum,
$\sum_{w\in A}\tau_w$,
can be close to 1 as much as we like (taking for instance
$l\to\infty$).
As we noticed above,
$\De_w>0$
unless
$s_w=t_w$.
So we always have
$\la(F,l,\al)>1$.
The point is that
$\la(F,l,\al)$
is separated from 1 uniformly over all
$l\ge l_0$, $\al\ge\al_0$.

Consider
$A_0\sub A$
consisting of all
$w\in A$
with
$|t_w-o_0|\ge 1$.
Then
$\De_w\ge\de_0>0$
for all
$w\in A_0$,
where
$\de_0$
depends only on
$l_0$, $\al_0$.
We claim that
\begin{equation}\label{equ:1}
\sum_{w\in A_0}\tau_w\ge m_0>0,
\end{equation}
where
$m_0=m_0(\wt F,l_0)$
is independent of
$l$.

Assuming (\ref{equ:1}), we have
\begin{eqnarray*}
\la(F,l,\al)&=&\sum_{w\in A}e^{\De_w}\tau_w
  \ge\sum_{w\in A_0}e^{\De_w}\tau_w+
     \sum_{w\in A\sm A_0}\tau_w\\
 &\ge& e^{\de_0}\sum_{w\in A_0}\tau_w+
      \sum_{w\in A\sm A_0}\tau_w\\
 &=&(e^{\de_0}-1)\sum_{w\in A_0}\tau_w+
   \sum_{w\in A}\tau_w\ge(e^{\de_0}-1)m_0+1=:\la_0>1.\\
\end{eqnarray*}

It remains to prove (\ref{equ:1}). Let
$o_w'\in w$
be the point closest to
$o_0\in w_0$.
Then
$\wt l_w\le|o-o_w'|\le l+|o_0-o_w'|$,
thus
$\wt l_w-l\le\dist(o_0,w)$. Hence
$\tau_w\ge e^{-\dist(o_0,w)}\ge\frac{\wt\psi_w}{4}$,
where
$w$
is observed from
$o_0$
under the angle
$\wt\psi_w$. 
The Lebesgue measure class on
$\dinf\hyp^2$
is independent of the choice of origin, thus
$\sum_{w\in A}\wt\psi_w=\pi$. 

Since
$l\ge l_0$,
for a sufficiently small
$m_0=m_0(S,l_0)>0$
the sectors
$S^+(m_0)$, $S^-(m_0)$,
defined below, intersect no
$w\in A\sm A_0$.
Here is the definition of
$S^{\pm}(m_0)$.
The common vertex
$o_0$
of
$S^{\pm}(m_0)$
divides the line
$w_0$
into two opposite rays
$w_0^{\pm}$.
The sectors
$S^{\pm}(m_0)\sub\hyp^2$
are bounded by the rays
$w_0^{\pm}$, $s^{\pm}(m_0)$,
where
$\angle_{o_0}(s^{\pm}(m_0),w_0^{\pm})=2m_0$,
and
$s^{\pm}(m_0)\cap F\not=\es$.

Therefore, it follows from
$\sum_{w\in A}\wt\psi_w=\pi$
that
$\sum_{w\in A_0}\tau_w\ge m_0$,
which completes the proof.
\end{proof}

\section{Accumulating procedure}\label{sect:accpro}

To describe the accumulating procedure we need some information 
about metric structure of the universal covering
$X$
of
$(M,g)$,
where
$g$
is a geometrization.

\subsection{Metric structure of the universal covering}

Recall (see, for example, \cite{BS}, \cite{CK}) that 
$X$
can be represented as the countable union
$X=\cup_vX_v$
of blocks, where each
$X_v$
is a closed convex subset in 
$X$
isometric to the metric product
$F_v\times\R$,
and 
$F_v$
is the universal covering of a compact hyperbolic surface
$S_v$
with geodesic boundary. Every two blocks are either disjoint
or intersect over a boundary component which is a 2-flat in
$X$
separating them and consequently no three blocks have a point
in common. The 2-flats in
$X$,
which separate blocks, are called the {\it walls}. A wall
$w$
common for blocks
$X_v$, $X_{v'}$
covers a 2-torus
$e\sub M$,
which separates (may be locally) the blocks
$M_v=\pi(X_v)$, $M_{v'}=\pi(X_{v'})$
of
$M$.
The metric decompositions
$X_v=F_v\times\R$, $X_{v'}=F_{v'}\times\R$
do not agree on
$w$,
and their 
$\R$-factors
induce two fibrations of
$w$
by parallel geodesics. We denote by
$\al_w$
the angle between these fibrations,
$0<\al_w\le\pi/2$.
Since
$M$
is compact and the set 
$E$
of separating tori in
$M$
is finite, we have
$\al_0:=\inf_w\al_w>0$,
where the infimum is taken over all the walls in
$X$.

\subsection{Modified Poincar\'e series}

We fix a wall
$w^*\sub X$,
take a block
$X_{v^*}\sub X$,
for which
$w^*$
is a boundary wall, and take a base point
$x_0\in w^*$.
We denote by
$W_0$
the set of the boundary walls of
$X_{v^*}$
different from
$w^*$;
by
$W_n$, $n\ge 1$
the set of walls in
$X$
at the combinatorial distance
$n+1$
from
$w^*$,
that is,
$w\in W_n$
if and only if any geodesic segment in
$X$
between
$x_0$
and
$w$
intersets
$n$
walls over its interior, including a wall from
$W_0$.
Note that
$W=\cup_{n\ge 0}W_n$
consists of all walls in
$X$,
which lie on one and the same side of
$w^*$
as
$X_{v^*}$.
Now, we define a modified Poincar\'e series
as follows
$$\cP_W(t)=\sum_{w\in W}e^{-t\dist(x_0,w)}.$$
Comparing
$P_W$
with
$\cP(t)=\sum_{\ga\in\Ga}e^{-t|x_0-\ga x_0|}$,
one easily obtains from triangle inequality that
$\cP(t)\ge e^{-D}\cP_W(t)$,
where
$D>0$
is the maximal diameter of the tori
$e\in E$.
Recall that the critical exponent of
$\cP$
is defined as the infimum of
$t\in\R$
for which
$\cP(t)<\infty$.
Therefore, the critical exponent 
$\ov h$
of
$\cP_W$
satisfies
$\ov h\le h(g)$,
and to prove Theorem~\ref{Thm:main} it suffices to show that
$\ov h>1$.

\subsection{Special broken-geodesic paths}

\begin{Pro}\label{Pro:brgeo} For each
$n\ge 0$
we have
$$\cP_n(1):=\sum_{w\in W_n}e^{-\dist(x_0,w)}\ge\frac{\pi}{4}\la_0^n,$$
where
$\la_0>1$
is the constant from Lemma~\ref{Lem:locest}.
\end{Pro}

\begin{proof} Induction over
$n$.
For each
$n\ge 0$
and each wall
$w\in W_n$
we produce a broken geodesic
$\xi_w$
in
$X$
between
$x_0$
and
$w$
as follows. For
$w\in W_0$
we put
$\xi_w=x_0x_w$,
where
$x_w\in w$
is the point closest to
$x_0$.
Note that
$x_0x_w$
lies in a horizontal slice
$F_{v^*}\times\{r_0\}$
of the block
$X_{v^*}=F_{v^*}\times\R$.

Assume that for all
$k$, $0\le k\le n-1$
and all
$w\in W_k$
the broken geodesic 
$\xi_w$
is already defined,  
$\xi_w$
is a juxtaposition
$\eta_{w_0}\eta_{w_1}\dots\eta_{w_k}'$
of geodesic segments with the sequence
$w_i\in W_i$
of walls leading to
$w=w_k$.
Moreover, we assume that each segment of
$\xi_w$
lies in a block and connects its different boundary components, 
the last segment
$\eta_{w_k}'$
lies in a horizontal slice of the block, which contains it,
and
$\eta_{w_k}'$
is orthogonal to the wall
$w$.

Take
$w\in W_n$.
Then there is a unique
$\ov w\in W_{n-1}$
which precedes
$w$.
By the assumption, the last edge
$\eta_{\ov w}'=s_{\ov w}x_{\ov w}$
of
$\xi_{\ov w}$
lies in a horizontal slice 
$F_{v'}\times\{r_{v'}\}$
of its block 
$X_{v'}=F_{v'}\times\R$
and it is
orthogonal to the wall
$\ov w$
(at the end point
$x_{\ov w}$).
The important feature of our construction is that 
$\xi_w$
contains all segments of
$\xi_{\ov w}$
but the last one
$\eta_{\ov w}'$.

Let
$X_v\sub X$
be the other block adjacent to
$\ov w$,
in particular,
the walls
$w$, $\ov w$
are its boundary components. Recall that the metric 
splitting
$X_v=F_v\times\R$
does not agree with that of
$X_{v'}$
along of
$\ov w$.
Let
$F_v\times\{r_v\}$
be the horizontal slice of
$X_v$
which contains
$x_{\ov w}$
on the corresponding boundary component. The boundary lines
of the slices
$F_{v'}\times\{r_{v'}\}$, $F_v\times\{r_v\}$, which lie in
$\ov w$,
contain
$x_{\ov w}$
and form an angle
$\al_{\ov w}\in(\al_0,\pi/2]$.

To construct
$\xi_w$
we do the following. We take an isometric copy
$F_v\sub\hyp^2\times\{r_{v'}\}$
of
$F_v\times\{r_v\}$
(by rotating the last by the angle
$\al_{\ov w}$),
where
$F_{v'}\times\{r_{v'}\}\sub\hyp^2\times\{r_{v'}\}$,
so that
$F_{v'}\times\{r_{v'}\}$
and
$F_v$
are sitting in the hyperbolic plane
$\hyp^2\times\{r_{v'}\}$
and they are adjacent along the common boundary component,
for which we use the same notation
$\ov w$.
Now, we connect the initial point
$s_{\ov w}$
of the last edge
$\eta_{\ov w}'\sub\xi_{\ov w}$
with the boundary component of
$F_v$
corresponding to
$w$
by the shortest geodesic segment
$s_{\ov w}x_w'\sub\hyp^2\times\{r_{v'}\}$
and take
$t_w=s_{\ov w}x_w'\cap\ov w$.
The segment
$t_wx_w'$
turned back to
$F_v\times\{r_v\}$
gives the last segment
$\eta_w'=s_wx_w$
of
$\xi_w$.
So
$s_w\in\ov w\cap F_v\times\{r_v\}$, 
$|s_w-x_{\ov w}|=|t_w-x_{\ov w}|$
and
$\eta_w'\sub F_v\times\{r_v\}$
is orthogonal to
$w$
at
$x_w$.
To complete the construction of
$\xi_w$,
we delete the last segment
$\eta_{\ov w}'\sub\xi_{\ov w}$
replacing it by
$\eta_{\ov w}\eta_w'$,
where
$\eta_{\ov w}=s_{\ov w}s_w\sub X_{v'}$.
Clearly, so constructed
$\xi_w$
has all properties advertised above.

Let
$l=|s_{\ov w}-x_{\ov w}|=L(\eta_{\ov w}')$
be the length of the last segment of
$\xi_{\ov w}$,
$l_w'=L(\eta_{\ov w})$, $l_w''=L(\eta_w')$.
Then we have
$$L(\xi_w)=L(\xi_{\ov w})-l+l_w'+l_w''=
   L(\xi_{\ov w})-l+\wt l_w-\De_w,$$
where
$\wt l_w=L(s_{\ov w}x_w')$,
$\wt l_w'=L(s_{\ov w}t_w)$
so that
$\wt l_w=\wt l_w'+l_w''$,
and
$\De_w=\wt l_w'-l_w'$.
We obtained the same configuration which was studied  in
sect.~\ref{sect:locest}, and we are going to apply 
Lemma~\ref{Lem:locest} to estimate
$\cP_n(1)$
from below. Note that
$l$
being the length of a segment in
$X$
connecting different boundary components of a block is
separated from 0 by some positive constant
$l_0$,
which depends only on
$M$, $l\ge l_0>0$.

For
$n=0$
we have
$$\cP_0(1)=\sum_{w\in W_0}e^{-|x_0-x_w|}=
   \sum_{w\in W_0}\tan\frac{\psi_w}{4}\ge\frac{\pi}{4},$$
where
$\psi_w$
is the angle under which the boundary component
$w$
of
$X_{v^*}$
is observed from
$x_0$
(in the horizontal direction).

By inductive assumption we have
$$\cP_{n-1}(1)\ge\sum_{w\in W_{n-1}}e^{-L(\xi_w)}\ge\frac{\pi}{4}\la_0^{n-1}.$$
Represent
$W_n=\cup_{\ov w}W_{n,\ov w}$,
where the union is taken over all
$\ov w\in W_{n-1}$
and each
$w\in W_{n,\ov w}$
follows
$\ov w$.
Then
$$\cP_n(1)\ge\sum_{w\in W_n}e^{-L(\xi_w)}
  =\sum_{\ov w\in W_{n-1}}e^{-L(\xi_{\ov w})}
  e^l\sum_{w\in W_{n,\ov w}}e^{\De_w}e^{-\wt l_w}.$$
Applying Lemma~\ref{Lem:locest} with
$A=W_{n,\ov w}$,
we obtain
$$\cP_n(1)\ge\la_0\sum_{\ov w\in W_{n-1}}e^{-L(\xi_{\ov w})}
  \ge\frac{\pi}{4}\la_0^n,$$
which completes the proof.
\end{proof}

\section{Self-similarity argument}\label{sect:selfsim}

The constant
$\la_0>1$
from Proposition~\ref{Pro:brgeo} depends only on some metric data of
$M$.
Thus there is
$n\in\N$, $n=n(M)$,
so that
$\frac{\pi}{4}\la_0^n>1$.
It follows from Proposition~\ref{Pro:brgeo} that
$$\cP_n(\ov h)=\sum_{w\in W_n}e^{-\ov h\dist(x_0,w)}\ge 1$$
for some
$\ov h>1$.
Furthermore, taking 
$n$
sufficiently large, we can find
$\ov h>1$
with
$\cP_n(\ov h)\ge 1$
for any choice of the initial block
$X_{v^*}$,
its wall
$w^*$
and the base point
$x_0\in w^*$,
since the set of choices up to isometries of
$X$
is compact.

We fix
$n\in\N$
with this property and select a subset
$W^*\sub W$, $W^*=\cup_{k\ge 1}W_k^*$,
where
$W_k^*=W_{kn}$.
The set
$W_1^*$
serves as the generating set for
$W^*$.
Connecting
$x_0$
with each wall
$w\in W_1^*$
by the shortest geodesic segment
$x_0x_w$,
we obtain new base points
$x_w\in w$
(these
$x_w$
may be different from
$x_w$
constructed in the proof of Proposition~\ref{Pro:brgeo}).
By induction, we find a base point
$x_w\in w$
for each
$w\in W_k^*$, $k\ge 1$
with the property 
$\dist(x_{\ov w},w)=|x_{\ov w}-x_w|$,
where
$\ov w\in W_{k-1}^*$
precedes
$w$.
Furthermore, by the choice of
$n$
and
$\ov h$
we have
$$\sum_we^{-\ov h|x_{\ov w}-x_w|}\ge 1,$$
for each
$\ov w\in W_{k-1}^*$,
where the summation is taken over all
$w\in W_k^*$
which follow
$\ov w$.
Since
$\dist(x_0,w)\le\dist(x_0,x_w)\le\dist(x_0,x_{\ov w})+|x_{\ov w}-x_w|$,
we obtain
\begin{eqnarray*}
 \cP_{W_k^*}(\ov h)&=&\sum_{w\in W_k^*}e^{-\ov h\dist(x_0,w)}\\
 &\ge&\sum_{w\in W_{k-1}^*}e^{-\ov h\dist(x_0,x_w)}
  \ge\sum_{w\in W_n}e^{-\ov h\dist(x_0,x_w)}\ge 1\\
\end{eqnarray*}
for each
$k\ge 1$.
Therefore, the modified Poincar\'e series
$$\cP_W(\ov h)\ge\cP_{W^*}(\ov h)=\sum_{k\ge 1}\cP_{W_k^*}(\ov h)$$
diverges at
$\ov h$
and hence
$h(g)\ge\ov h>1$.
This completes the proof of Theorem~\ref{Thm:main}.

\end{document}